\documentclass{amsart}

\usepackage{amssymb,amsmath,latexsym,euscript,epsfig}
\usepackage{graphics}

 \newtheorem{theorem}{Theorem}
 \newtheorem{conjecture}[theorem]{Conjecture}
 \newtheorem{lemma}[theorem]{Lemma}
 
 \theoremstyle{definition}

\newtheorem*{remark}{Remark}

\def\A{{\mathcal A}}
\def\AA{{\mathfrak A}}

\def\FF{{\mathcal F}}

\def\L{{\mathcal L}}

\def\P{{\mathcal P}}

\def\T{{\mathbb T}}

\def\Z{{\mathbb Z}}

\def\ls{\le}

\def\w{{\bf w}}
\def\wB{{\widetilde{B}}}
\def\wx{{\widetilde{\bf x}}}
\def\x{{\bf x}}
\def\y{{\bf y}}

\def\cf{{\operatorname{cf}}}

\def\diag{\operatorname{diag}}

\def\pr{{\operatorname{pr}}}

\def\:{{:\ }}

\begin{document}
\title{On the properties of the exchange graph of a cluster algebra}

\author{Mikhail Gekhtman}
\address{Department of Mathematics, University of Notre Dame, Notre Dame,
Indiana 46556}
\curraddr{Department of Mathematics, University of Notre Dame, Notre Dame,
Indiana 46556}
\email{mgekhtma@nd.edu}
\thanks{The authors were supported in part by BSF Grant \#2002375.}
\thanks{The first author was supported in part by NSF Grant \#0400484.}

\author{Mikhail Shapiro}
\address{Department of Mathematics, Michigan State University, East Lansing,
Michigan 48823}
\email{mshapiro@math.msu.edu}
\thanks{The second author was supported in part by PHY Grant \#0555346 and DMS Grant \#0401178.}

\author{Alek Vainshtein}
\address{Department of Mathematics AND Department of Computer Science, University of Haifa, Haifa,
Mount Carmel 31905, Israel}
\email{alek@cs.haifa.ac.il}

\date\today

\begin{abstract}
We prove a conjecture about the vertices and edges of the exchange graph of a cluster
algebra $\A$ in two cases: when $\A$ is of geometric type and when $\A$ is arbitrary and
its exchange matrix is nondegenerate. In the second case we also prove that the exchange
graph does not depend on the coefficients of $\A$. Both conjectures were formulated
recently by Fomin and Zelevinsky.
\end{abstract}

\maketitle

\section{Main definitions and results}

A cluster algebra  is an axiomatically defined commutative ring
equipped with a distinguished set of generators (cluster
variables). These generators are subdivided into overlapping subsets (clusters)
of the same cardinality that are connected via sequences of birational transformations 
of a particular kind, called cluster
transformations. Transfomations of this kind can be observed in many areas of
mathematics (Pl\"ucker relations, Somos sequences and Hirota equations,
to name just a few examples).

Cluster algebras were initially
introduced in \cite{CAI} to study total positivity and (dual)
canonical bases in semisimple algebraic groups.
The rapid development of the cluster algebra
theory revealed relations between cluster algebras and Grassmannians,
quiver representations, generalized associahedra, Teichm\"uller
theory, Poisson geometry  and many other branches of mathematics, see
\cite{Ze} and references therein. In the present paper we prove some
conjectures on the general structure of  cluster algebras formulated
by~Fomin and~Zelevinsky in \cite{CAC}.

To state our results, we recall the definition of a cluster algebra; for
details see \cite{CAI,CAIV}.

Let $\P$ be a semifield, that is, a torsion-free multiplicative abelian group 
endowed with an additional operation
$\oplus$, which is commutative, associative and distributive with respect to the multiplication.
As an {\em ambient field\/} we take
the field $\FF$ of rational functions in $n$ independent variables with
coefficients in the field of fractions of the integer group ring
$\Z\P$. A square matrix $B$ is called {\em skew-symmetrizable\/} if $DB$ is skew-symmetric for a
positive diagonal matrix $D$.
A {\em seed\/} is a triple $\Sigma=(\x,\y,B)$, where $\x=(x_1,\dots,x_n)$ is a transcendence basis
 of $\FF$ over the field of
fractions of $\Z\P$, $\y=(y_1,\dots,y_n)$ is an $n$-tuple of elements of $\P$ and $B$ is a
skew-symmetrizable integer $n\times n$ matrix. The components $\x$, $\y$ and
$B$ of the seed are called 
the {\em cluster\/},
the {\em coefficient tuple\/} and the {\em exchange matrix\/}, respectively; the entries of $\x$ are
called {\em cluster variables}.

A {\em seed mutation\/} in direction $k\in [1,n]$ takes $\Sigma$ to an {\em adjacent\/} seed
$\Sigma'=(\x',\y',B')$ whose components are defined as follows. The adjacent cluster $\x'$ is
given by $\x'=(\x\setminus\{x_k\})\cup\{x'_k\}$,
where the new cluster variable $x'_k$ is defined by the {\em exchange relation}
\begin{equation}\label{newexch}
x_kx'_k=\frac{y_k}{y_k\oplus 1}\prod_{b_{ki}>0}x_i^{b_{ki}}+
\frac{1}{y_k\oplus 1}\prod_{b_{ki}<0}x_i^{-b_{ki}},
\end{equation}
where, as usual, the product over the empty set is assumed to be
equal to~$1$.
The adjacent coefficient tuple $\y'$ is given by
\begin{equation}\label{newcoefch}
y'_j=\begin{cases} y_k^{-1} &\text{if $j=k$},\\
                   y_jy_k^{b_{jk}}(y_k\oplus 1)^{-b_{jk}} &\text{if $j\ne k$ and $b_{jk}>0$},\\
                   y_j (y_k\oplus 1)^{-b_{jk}} &\text{if $j\ne k$ and $b_{jk}\ls 0$}.
     \end{cases}
\end{equation}
Finally, the adjacent exchange matrix $B'$ is given by
\begin{equation}\label{matmut}
b'_{ij}=\begin{cases}
         -b_{ij}, & \text{if $i=k$ or $j=k$;}\\
                 b_{ij}+\displaystyle\frac{|b_{ik}|b_{kj}+b_{ik}|b_{kj}|}2,
                                                  &\text{otherwise.}
        \end{cases}
\end{equation}

\begin{remark}
Note that our encoding of exponents in relations~\ref{newexch} and~\ref{newcoefch} by the entries
of $B$ differs from the one used in~\cite{CAI}-\cite{CAIV} by transposition of subscripts. 
In other words, our exchange matrix is the transpose of the one used in~\cite{CAI}-\cite{CAIV}.
\end{remark}

Two seeds are called {\em mutation equivalent\/} if they can be obtained one 
from another by a sequence of seed mutations.
The {\em cluster algebra\/} $\A=\A(\Sigma)$
associated with $\Sigma$ is the $\Z\P$-subalgebra of $\FF$
generated by all cluster variables in all seeds mutation
equivalent to $\Sigma$; $n$ is said to be the {\em rank\/} of $\A(\Sigma)$.

We say that $\A$ is of {\em geometric type\/} if $\P$ is a {\em tropical semifield\/}, 
that is, if $\oplus$ is defined by
\[
\prod_{i=1}^m g_i^{a_i}\oplus \prod_{i=1}^m g_i^{b_i}=
\prod_{i=1}^m g_i^{\min(a_i,b_i)}.
\]
For a cluster algebra of geometric type, each entry of the coefficient tuple $\y$ can be represented as
\begin{equation}\label{yfromg}
y_i=\prod_{i=1}^m g_j^{a_{ij}}
\end{equation}
for some integer $a_{ij}$. Re-denote the generators $g_1,\dots,g_m$ by $x_{n+1},\dots,x_{n+m}$  and call
them {\em stable variables\/}; together with the cluster variables, they 
form an {\em extended cluster\/}.
Define an {\em extended exchange matrix\/} $\wB$ as an $n\times (n+m)$ matrix
whose $(i,j)$th entry equals $b_{ij}$ for $i\in [1,n]$ and $a_{i-n,j}$ for $i\in [n+1,n+m]$.
For an arbitrary matrix $A$, denote by
$A[k;l]$ the submatrix of $A$ occupying the first $k$ rows and
the first $l$ columns; with this notation, the relation between $B$ and
$\wB$ is given by $B=\wB[n;n]$.

Transformation rules~(\ref{newexch}) can be rewritten  for cluster algebras 
of geometric type as
\begin{equation}\label{exchange}
x_kx'_k=\prod_{\substack{1\le i\le n+m\\  b_{ki}>0}}x_i^{b_{ki}}+
       \prod_{\substack{1\le i\le n+m\\  b_{ki}<0}}x_i^{-b_{ki}}.
\end{equation}
Transformation 
rules~(\ref{newcoefch}) imply that mutations of the extended exchange matrix are governed by the
same  matrix mutation rules~(\ref{matmut}). The corresponding cluster
algebra of geometric type is denoted $\A(\wB)$.

A convenient tool in dealing with cluster algebras is the $n$-regular tree $\T_n$.
Its vertices correspond to seeds, and two vertices are connected by an edge labeled
by $k$ if and only if the corresponding seeds are adjacent in direction $k$.
The edges of $\T_n$ are thus labeled by the numbers $1,\dots,n$ so that the $n$ edges
emanating from each vertex receive different labels.
Two seeds differing from each
other by an arbitrary permutation of $\x$
and the corresponding permutation of $\y$ and of the rows and columns of
$\wB$ are called {\it equivalent}. The {\em exchange
graph\/} of a cluster  algebra is defined as the quotient of the 
tree $\T_n$ modulo this equivalence relation.

 The following conjectures about the exchange graph of a cluster algrebra
were formulated in \cite{CAC}.

\begin{conjecture}\label{conj:Bdefgraph}
 The exchange graph of a cluster algebra depends only on the initial
exchange matrix $B$.
\end{conjecture}

\begin{conjecture}\label{conj:cldefseed}
 Every seed is uniquely defined by its cluster; thus, the vertices of the
exchange graph can be
identified with the clusters, up to a permutation of cluster variables.
\end{conjecture}

\begin{conjecture}\label{conj:adjclusters}
 Two clusters are adjacent in the exchange graph if and only if they have
exactly $n-1$ common cluster variables.
\end{conjecture}

All three conjectures were proved for cluster algebras
of finite type in  \cite{CAII}, and for cluster algebras
associated with triangulations of two-dimensional surfaces in \cite{FST}.
Besides, in \cite{BMRT}
 Conjectures~\ref{conj:cldefseed} and~\ref{conj:adjclusters} were proved
for cluster algebras of
 geometric type with no stable variables with a skew-symmetric exchange
matrix satisfying an additional acyclicity condition.

 In this note we prove the following three results. Denote by $\AA(B)$ the
family of all cluster algebras
 obtained by fixing the initial exchange matrix $B$ and varying the
semifield $\P$, the choice of the initial coefficient $n$-tuple $y$,
and the choice of the initial cluster $x$ (i.e.,  the choice of a
transcendence basis of $\FF$ over the field of fractions of $\Z\P$).


\begin{theorem}\label{thm:cldefseed}
Let a cluster algebra $\A\in\AA(B)$ satisfy one of the following two
conditions:

(i) $\A$ is of geometric type;

(ii) $\A$ is arbitrary and $B$ is nondegenerate.

Then every seed in $\A$ is uniquely defined by its cluster
{\rm(}in other words, Conjecture~\ref{conj:cldefseed} holds true for
$\A${\rm)}.
\end{theorem}

\begin{theorem}\label{thm:adjclusters}
Suppose that every seed in a cluster algebra $\A$ is uniquely defined by
its cluster. Then two clusters are adjacent in the exchange graph of $\A$
if and only if they have exactly $n-1$ common variables. In particular,
Conjecture~\ref{conj:adjclusters} holds true if $\A$ satisfies one of
the conditions {\rm(i),} {\rm(ii)} of Theorem~\ref{thm:cldefseed}.
\end{theorem}

\begin{remark} Note that we do not claim that the whole
cluster complex is determined by set-theoretic combinatorics of
subsets (clusters) of the ground set of all cluster variables,
although we strongly believe that it is true.
Theorem~\ref{thm:adjclusters}
states only that the exchange graph of a
cluster algebra is determined by this combinatorics. In other words, we
can restore 
which maximal simplices are adjacent by a codimension one face.
\end{remark}

\begin{theorem}\label{thm:Bdefgraph}
 Let $B$ be nondegenerate, then the exchange graphs of all cluster
algebras in $\AA(B)$ coincide {\rm(}in other words,
Conjecture~\ref{conj:Bdefgraph} holds true{\rm)}.
\end{theorem}

We prove Theorems~\ref{thm:cldefseed}--\ref{thm:Bdefgraph} 
using the concept of a 2-form compatible with the cluster
algebra structure \cite{GSV1,GSV2}. Under certain nondegeneracy conditions, this approach
has already proved useful in realizing cluster algebras as coordinate rings
of rational Poisson manifolds.

The work on this
project started during the visit of  M.~G.~and M.~S.~to the Haifa
Interdisciplinary Research Center for Advanced Computer Science
in Fall 2006 and completed during the visit of A.~V.~to the Michigan State University in Spring 2007.
The authors are grateful to these institutions for hospitality and stimulating atmosphere.
The authors also thank
A.~Zelevinsky for attracting their attention to \cite{BMRT}.

\section{Proof of Theorem~\ref{thm:cldefseed}}
In what follows we write $x_{i;\Sigma}$, $y_{i;\Sigma}$ and $b_{ij;\Sigma}$ to indicate that
the corresponding cluster variables, coefficients and entries of the
exchange matrix are related to a seed $\Sigma$.

Observe that Conjecture~\ref{conj:cldefseed} is equivalent to the
following one:
suppose that the seeds 
$\Sigma_1$ and $\Sigma_2$ satisfy relations
\begin{equation}\label{samecl}
x_{i;\Sigma_2}=x_{\sigma(i);\Sigma_1}
\end{equation}
for some $\sigma\in S_n$ and any $i\in [1,n]$, then
\begin{equation}\label{samey}
y_{i;\Sigma_2}=y_{\sigma(i);\Sigma_1}
\end{equation}
for $i\in [1,n]$ and
\begin{equation}\label{sameb}
b_{ij;\Sigma_2}=b_{\sigma(i)\sigma(j);\Sigma_1}
\end{equation}
for any $i,j\in [1,n]$.

First, consider case (i) in Theorem~\ref{thm:cldefseed}.
Following \cite{GSV2}, we say that a closed rational differential
2-form $\omega$ 
on the $(n+m)$-dimensional
affine space is {\it compatible\/} with the cluster algebra $\A(\wB)$ if for any extended cluster
$\wx=(x_1,\dots,x_{n+m})$ one has
$$
\omega=\sum_{i,j=1}^{n+m}\omega_{ij}\frac{d x_i}{x_i}\wedge \frac{d x_j}{x_j},
$$
where $\omega_{ij}$ are constants (recall that $x_i$'s are rational
functions in the initial cluster variables).
The matrix $\Omega=(\omega_{ij})$ is called the {\it
coefficient matrix\/} of $\omega$ (with respect to $\wx$); evidently,
$\Omega$ is skew-symmetric.

A square matrix $A$ is {\it decomposable\/} if there
exists a permutation matrix $P$ such that $PAP^T$ is a block-diagonal
matrix, and {\it indecomposable\/} otherwise;
$\rho(A)$ is defined as the maximal number of diagonal blocks in $PAP^T$.
The partition into blocks defines an obvious equivalence relation
$\sim$ on the rows (or columns) of $A$.

The following result is a generalization of Theorem~2.1 of \cite{GSV2}.

\begin{theorem}\label{thm:comp2form}
Assume that $\wB$ is $D$-skew-symmetrizable and does not have zero rows.
Then all rational closed 2-forms compatible with $\A(\wB)$ form a vector
space of dimension $\rho(B)+\binom{m}2$, where $B=\wB[n;n]$.
Moreover, the coefficient matrices of these 2-forms with respect to $\wx$
are characterized by the equation
$\Omega[n;n+m]=\Lambda D\wB$, where
$\Lambda=\diag(\lambda_1,\dots,\lambda_n)$ with
$\lambda_i=\lambda_j$ whenever $i\sim j$. In particular, if $B$ is
indecomposable, then $\Omega[n;n+m]=\lambda D\wB$.
\end{theorem}

\proof
Indeed, let $\omega$ be a 2-form compatible with $\A(\wB)$. Then
$$
\omega=\sum_{j,k=1}^{n+m}\omega_{jk}\frac{dx_j}{x_j}\wedge \frac{d x_k}{x_k}
=\sum_{j,k=1}^{n+m}\omega'_{jk}
\frac{dx'_j}{x'_j}\wedge \frac{d x'_k}{x'_k},
$$
where $x'_j$ is given by~(\ref{exchange})
and $\omega'_{jk}$ are the coefficients of
$\omega$ with respect to $\wx'$. Recall that the only variable in the
extended cluster
$\wx'$ different from the corresponding variable in $\wx$ is
$x_i$, and
\begin{equation*}
\frac{dx'_i}{x'_i}=-\frac{dx_i}{x_i} +
\sum_{b_{ik}>0}\frac{b_{ik}}{1+\prod_{k=1}^{n+m}x_k^{-b_{ik}}}\frac{dx_k}{x_k}
-\sum_{b_{ik}<0}\frac{b_{ik}}{1+\prod_{k=1}^{n+m}x_k^{b_{ik}}}\frac{dx_k}{x_k}.
\end{equation*}

Thus, for any $j\in [1,n+m]$ we immediately get
\begin{equation}\label{eq:2.3a}
\omega'_{ij}=-\omega_{ij}.
\end{equation}

Next, consider any pair $j,k\ne i$ such that both $b_{ij}$ and $b_{ik}$
are nonnegative, and at least one of the two is positive. Then
$$
\omega_{jk}=\omega'_{jk}+
\frac{\omega'_{ik}b_{ij}+\omega'_{ji}b_{ik}}
{1+\prod_{k=1}^{n+m}x_k^{-b_{ik}}}.
$$
This equality can only hold if
$\omega'_{ik}b_{ij}+\omega'_{ji}b_{ik}=0$, which by~(\ref{eq:2.3a}) is
equivalent to $\omega_{ij}b_{ik}=\omega_{ik}b_{ij}$.
If both $b_{ij}$ and $b_{ik}$ are positive, this
gives
\begin{equation}\label{eq:2.4}
\frac{\omega_{ij}}{b_{ij}}=\frac{\omega_{ik}}{b_{ik}}=\mu_{i}.
\end{equation}
Otherwise, if, say, $b_{ij}=0$, one gets $\omega_{ij}=0$.
Besides, in any case $\omega'_{jk}=\omega_{jk}$.

Similarly, if both $b_{ij}$ and $b_{ik}$ are nonpositive,
and at least one of the two is negative, then
$$
\omega_{jk}=\omega'_{jk}-
\frac{\omega'_{ik}+\omega'_{ji}b_{ik}}
{1+\prod_{k=1}^{n+m}x_k^{b_{ik}}},
$$
and hence the same relations as above hold true.

Finally, let
$b_{ij}\cdot b_{ik}<0$,  say, $b_{ij}>0$ and
$b_{ik}<0$; then
$$
\omega_{jk}=\omega'_{jk}+
\frac{\omega'_{ik}b_{ij}\prod_{k=1}^{n+m}x_k^{b_{ik}}-\omega'_{ji}b_{ik}}
{1+\prod_{k=1}^{n+m}x_k^{b_{ik}}},
$$
which again leads to~(\ref{eq:2.4});
the only difference is that in this case
$\omega'_{jk}=\omega_{jk}+\omega_{ik}b_{ij}$.

We have thus obtained that
 $\Omega[n;n+m]=\diag(\mu_1,\dots,\mu_n)\wB$. Recall that
 $\Omega$ is skew-symmetric; however, any skew-symmetrizer for $\wB$
can be written as $\Lambda D$, where
$\Lambda=\diag(\lambda_1,\dots,\lambda_n)$ and
$\lambda_i=\lambda_j$ whenever $i\sim j$, which completes the proof.
\endproof

Assume now that the initial extended matrix $\wB$ does not have zero rows.
Then by Theorem~\ref{thm:comp2form} we can pick a rational closed 2-form
$\omega$ compatible with $\A$ so that for any seed $\Sigma$
holds
\begin{equation}\label{comp2form}
B_\Sigma=D^{-1}\Omega_\Sigma[n;n+m],
\end{equation}
where $D$ is a skew-symmetrizer of $B$ and $\Omega_\Sigma$ stands for the
coefficient matrix of $\Omega$ with respect to $\wx_\Sigma$. Since~(\ref{samecl})
implies $\omega_{ij;\Sigma_2}=\omega_{\sigma(i)\sigma(j);\Sigma_1}$, we get
\begin{equation}\label{adh71}
b_{\sigma(i)\sigma(j);\Sigma_1}=d^{-1}_{\sigma(i)}\omega_{\sigma(i)\sigma(j);\Sigma_1}=
d^{-1}_{\sigma(i)}\omega_{ij;\Sigma_2}=d_id^{-1}_{\sigma(i)}b_{ij;\Sigma_2}.
\end{equation}
Let us show that $d_i=d_{\sigma(i)}$. Put $\alpha=d_id_{\sigma(i)}^{-1}$; clearly, $\alpha>0$, since both $d_i$ and $d_{\sigma(i)}$
are positive. Let $\Sigma_2'$ be the seed adjacent to $\Sigma_2$ in direction $i$, and $\Sigma_1'$ be the seed adjacent to $\Sigma_1$
in direction $\sigma(i)$. Assume that $x_{i;\Sigma_2'}x_{i;\Sigma_2}=\Pi_1+\Pi_2$, where $\Pi_1,\Pi_2$ are monomials in the
variables of $\{x_{j;\Sigma_2}\: j\ne i\}$. Then~\eqref{samecl} and~\eqref{adh71} imply $x_{\sigma(i);\Sigma_1'}x_{\sigma(i);\Sigma_1}=
\Pi_1^\alpha+\Pi_2^\alpha$. Using~\eqref{samecl} once again we get
\begin{equation}\label{misha}
x_{\sigma(i);\Sigma_1'}=\frac{\Pi_1^\alpha+\Pi_2^\alpha}{\Pi_1+\Pi_2}x_{i;\Sigma_2'}.
\end{equation}

By the Laurent phenomenon \cite{CAI}, $x_{\sigma(i);\Sigma_1'}$ should be represented as a Laurent polynomial in 
$\{x_{j;\Sigma_2}\: j\ne i\}\cup \{x_{i;\Sigma_2'}\}$. By~\eqref{misha}, this is impossible for $0<\alpha<1$. Similarly,
$x_{i;\Sigma_2'}$ should be represented as a Laurent polynomial in $\{x_{\sigma(j);\Sigma_1}\: j\ne i\}\cup \{x_{\sigma(i);\Sigma_1'}\}$,
which is impossible for $\alpha>1$. Therefore, $\alpha=1$.

We thus proved that~(\ref {sameb}) holds true
for any $i,j\in [1,n]$. Similarly, taking into account~(\ref{comp2form}) and
$d_i=d_{\sigma(i)}$ we get $b_{i,n+j;\Sigma_2}=b_{\sigma(i),n+j;\Sigma_1}$
for any $i,j\in [1,n]$.
It remains to recall that for cluster algebras of geometric type the
$n$-tuple $\y_\Sigma$ is completely defined by
the matrix $B_\Sigma$ via~(\ref{yfromg}).

Assume now that $\wB$ has $k$ zero rows. In this case $\A$ is a direct product of $k$
copies of the cluster algebra of rank~1  and a cluster algebra defined by a
$(n-k)\times(n+m)$-submatrix of $\wB$ with no zero rows, for which the above reasoning applies.

Consider now case (ii). Define the {\em coefficient-free\/} cluster algebra
$\A_\cf(B)\in\AA(B)$ as the cluster algebra over the one-element semifield $\{1\}$ obtained
by mutations of the seed $(\x,1,B)$. Clearly, the map that takes all elements of $\P$
to~$1$ commutes with the exchange relation~(\ref{newexch}).
Therefore,~(\ref{samecl}) for an algebra $\A\in\AA(B)$ implies the same relation for
$\A_\cf(B)$. Since the latter is of geometric type, the result of case (i) applies,
and hence relation~(\ref{sameb}) holds true.

Following \cite{CAIV}, for an arbitrary cluster algebra $\A\in\AA(B)$ introduce the
$n$-tuple $\widehat{\y}_\Sigma=(\widehat{y}_{1;\Sigma},\dots,\widehat{y}_{n;\Sigma})$ by
\[
\widehat{y}_{j;\Sigma}=y_{j;\Sigma}\prod_{k=1}^nx_{k;\Sigma}^{b_{jk;\Sigma}}.
\]
For cluster algebras of geometric type, in particular, for the coefficient-free algebra
$\A_\cf(B)$, this $n$-tuple coincides with $\tau$-coordinates introduced
in \cite{GSV1}.
Further, define $Y_{i;\Sigma}(y_1,\dots,y_n)$ as a rational
subtraction-free expression for the coefficient
$y_{i;\Sigma}$ via the initial coefficients $y_1,\dots,y_n$. The evaluation of
$Y_{i;\Sigma}$ over $\P$ denoted by $Y_{i;\Sigma}|\P$ gives the value of $y_{i;\Sigma}$
in the cluster algerba over $\P$; this operation is well-defined since any subtraction-free identity
in the field of rational functions remains valid in any semifield (see~\cite{BFZ}, Lemma~2.1.6).

By Proposition~3.9 in \cite{CAIV},
\[
\widehat{y}_{j;\Sigma}=Y_{j;\Sigma}|_\FF(\widehat{y}_1,\dots,\widehat{y}_n).
\]
Therefore,~(\ref{samecl}) and~(\ref{sameb}) for $\A_\cf(B)$ imply
$Y_{i;\Sigma_2}(\tau_1,\dots,\tau_n)=Y_{\sigma(i);\Sigma_1}(\tau_1,\dots,\tau_n)$.

 Since $B$ is nondegenerate, Lemma~1.1 in \cite{GSV1} implies that the transformation
$\x\mapsto\tau$ is a bijection, and hence
$Y_{i;\Sigma_2}(z_1,\dots,z_n)=Y_{\sigma(i);\Sigma_1}(z_1,\dots,z_n)$ for any set of variables
$z_1,\dots,z_n$. Therefore,~(\ref{samey}) holds for any cluster algebra $\A\in\AA(B)$.

\begin{remark} An similar approach based on utilizing the bijection 
$\x\mapsto\tau$ is used in the proof of Proposition~2.7 in
\cite{FoG}.
\end{remark}

\section{Proof of Theorem~\ref{thm:adjclusters}}

Denote by $x_1,\dots,x_{n-1}$ the common variables in the two clusters and
by $\dot x$ and $\ddot x$ the remaining variables;
the clusters themselves will be denoted $\dot\x$ and $\ddot\x$, respectively.
By Theorem~3.1 of \cite{CAI}, $\ddot x$ can be written as a Laurent polynomial in
$x_1,\dots,x_{n-1},\dot x$. Since each cluster transformation is birational,
$\dot x$ enters this polynomial with exponent~$1$ or~$-1$; we write this as
$\ddot x=\L_0+L_1\dot x^{\pm1}$ where $L_0$ and $L_1$ are Laurent polynomials in
$x_1,\dots,x_{n-1}$. Denote by $\ddot x'$ the cluster variable in the
cluster 
adjacent to $\ddot\x$ that replaces $\ddot x$, then
\begin{equation}\label{xtodd}
\ddot x'=\frac{M+N}{L_0+L_1\dot x^{\pm1}},
\end{equation}
where $M$ and $N$ are monomials in $x_1,\dots,x_{n-1}$. Since $\ddot
x'$ is a Laurent polynomial in $x_1,\dots,x_{n-1},x$, we immediately get $L_0=0$.

We have to consider two cases: $\ddot x=L_1\dot x$ and $\ddot
x=L_1/\dot x$. 
In the first case we write down $\dot x=L_1^{-1}\ddot x$ and apply
once again  
Theorem~3.1 of \cite{CAI} to see that $L_1$ is a Laurent monomial
$M_+$. In the second case we use~(\ref{xtodd})
to get $\ddot x'=M_-\dot x$ for some Laurent monomial $M_-$. It
remains to prove that both 
$M_+$ and $M_-$ are identically equal to~$1$. In the first case this
would mean $\ddot\x$ and 
$\dot\x$ have $n$ common cluster variables, a contradiction. In the
second case this would mean that 
$\dot \x$ is adjacent to $\ddot\x$ as required.

Assume that a variable $x_i$ enters $M_+$ with a negative exponent
$-k$, $k>0$. Let $x'_i$ denote the 
cluster variable that replaces $x_i$ in the cluster adjacent to
$\dot\x$ in the corresponding 
direction. Then $x_ix_i'=M_i+N_i$, where $M_i$ and $N_i$ are monomials
in $x_1,\dots,x_{n-1},\dot x$. Clearly, at least one of $M_i$ and $N_i$ is 
nontrivial, since otherwise the cluster algebra of rank~1 generated by $x_i$ 
would split off, in a contradiction to the assumption that $x_i$ enters $M_+$.
Therefore,
\[
\ddot x=\frac{\dot x(x_i')^kM'_+}{(M_i+N_i)^k},
\]
where $M'_+$ is a Laurent monomial and the denominator  is a
nontrivial polynomial, 
which contradicts the Laurent property. To handle the case of a
positive exponent $k>0$ we write 
$\dot x=M_+^{-1}\ddot x$ and proceed in a similar way using the
variable $x_i''$  that 
replaces $x_i$ in the cluster adjacent to $\ddot\x$ in the
corresponding direction. 
We thus obtained that $M_+$ is a constant. It is an easy exercise to
prove that 
the constant is equal to 1. The case of $M_-$ is handled similarly.

\section{Proof of Theorem~\ref{thm:Bdefgraph}}
Let $B_\pr$ be the $n\times 2n$ matrix whose principal part equals $B$ and
the remaining part is the $n\times n$ identity matrix. The corresponding algebra of geometric
type $\A(B_\pr)\in\AA(B)$ is called the algebra with {\it principal coefficients\/} and is denoted by
$\A_\pr(B)$. Note that for this algebra initial coefficients 
$y_1,\dots,y_n$ coincide with the stable variables $x_{n+1},\dots,x_{2n}$.

The exchange graph of $\A'$ {\it covers\/} the exchange graph of $\A$ if the equivalence of
two seeds in $\A'$ implies the equivalence of the corresponding two seeds in $\A$.
It follows immediately from the definition that the exchange graph of any cluster algebra in $\AA(B)$
covers the exchange graph of the coefficient-free cluster algebra
$\A_\cf(B)$. By Theorem~4.6 in \cite{CAIV}, the exchange graph of 
$\A_\pr(B)$ covers the exchange graph of any $\A\in\AA(B)$. Therefore,
it suffices to prove that the exchange graphs for $\A_\pr(B)$ and $\A_\cf(B)$ coincide.

Define $X_{i;\Sigma}(x_1,\dots,x_n;y_1,\dots,y_n)$ as a rational 
function expressing the cluster variable $x_{i;\Sigma}$ in
$\A_\pr(B)$; 
further, define rational functions $G_{i;\Sigma}$ by
$G_{i;\Sigma}(x_1,\dots,x_n)=X_{i;\Sigma}(x_1,\dots,x_n;1,\dots,1)$. 
By Theorem~3.7 in \cite{CAIV}, $G_{i;\Sigma}$ express the variables 
$x_{i;\Sigma}$ via the initial variables in the coefficient-free algebra $\A_\cf(B)$.

By Theorem~\ref{thm:cldefseed}, we have to prove the following implication: if
\begin{equation}\label{sameG}
G_{i;\Sigma_2}=G_{\sigma(i),\Sigma_1}
\end{equation}
for some $\sigma\in S_n$, some seeds $\Sigma_1,\Sigma_2$ and any $i\in [1,n]$, then
\begin{equation*}\label{sameX}
X_{i;\Sigma_2}=X_{\sigma(i);\Sigma_1}
\end{equation*}
for any $i\in [1,n]$. The proof is based on the following lemma.

\begin{lemma}\label{lem:GtoX}
For any seed $\Sigma$ there exist Laurent monomials $M_{i;\Sigma}$ such that
\begin{multline}\label{eq:GtoX}
X_{i;\Sigma}(x_1,\dots,x_n;y_1,\dots,y_n)\\
=M_{i;\Sigma}(z_1,\dots,z_n)
G_{i;\Sigma}(x_1M_1(z_1,\dots,z_n),\dots, x_nM_n(z_1,\dots,z_n)),\qquad i\in [1,n],
\end{multline}
where $M_i=M_{i;\Sigma_0}$,  $\Sigma_0=(\x,\y,B)$ is the initial cluster  and $z_i=y_i^{1/\det B}$.
\end{lemma}

\begin{proof}
Consider toric actions on $\A_\pr(B)$ similar to those
introduced in \cite{GSV1}. Given an $n$-tuple of integer weights
$$
\w_\Sigma^1=(w^1_{1;\Sigma},\dots,w^1_{2n;\Sigma}),\dots, \w^n=(w^n_{1;\Sigma},\dots,w^n_{2n;\Sigma})
$$
for each seed $\Sigma$, define a {\em local\/} toric action by
\[
(x_{1;\Sigma},\dots,x_{2n;\Sigma})\mapsto
 (t_1^{w^1_{1;\Sigma}}\dots t_n^{w^n_{1;\Sigma}}x_{1;\Sigma},\dots,
 t_1^{w^1_{n;\Sigma}}\dots t_n^{w^n_{n;\Sigma}}x_{n;\Sigma}).
\]
Two local toric actions at adjacent seeds are compatible if they
commute 
with the transformation given by~(\ref{exchange}). If all local toric
actions 
are compatible, they determine a {\em global\/} toric action on 
$\A(\wB)$ called the extension of each of the local actions.

A slight modification of Lemma~2.3 in \cite{GSV1} guarantees that if
$$
\w^1=(w^1_1,\dots,w^1_{2n}), \dots, \w^n=(w^n_1,\dots,w^n_{2n})
$$
satisfy $\wB(\w^j)^T=0$ for all $j\in [1,n]$ then the local toric action at $\Sigma_0$
defined by $\w^1,\dots,\w^n$ can be extended to a global toric action.

Define the weights $\w^i$ as follows: the first $n$ entries of 
$\w^i$ constitute the $i$th row of the matrix $B^{-1}$ multiplied by $\det B$, while its
last $n$ entries constitute the $i$th row of the 
$n\times n$-matrix $\diag(-\det B,\dots,-\det B)$. 
Then the compatibility condition gives
\begin{eqnarray*}
X_{i;\Sigma}(x_1t_1^{w^1_1}\dots t_n^{w^n_1},\dots,x_nt_1^{w^1_n}\dots t_n^{w^n_n};
y_1t_1^{-\det B},\dots,y_nt_n^{-\det B})\\
=N_{i;\Sigma}(t_1,\dots,t_n)X_{i;\Sigma}(x_1,\dots,x_n;y_1,\dots,y_n)
\end{eqnarray*}
for some Laurent monomials $N_{i;\Sigma}$. Define monomials $M_i$ by
$$
M_i(t_1,\dots,t_n)=\prod_{j=1}^nt_j^{w^i_j}.
$$
Relation~(\ref{eq:GtoX}) follows from the above condition with
$t_i=z_i$ 
and $M_{i;\Sigma}=N_{i;\Sigma}^{-1}$.
\end{proof}

From~(\ref{sameG}) and~(\ref{eq:GtoX}) we get
$$
X_{i;\Sigma_2}=M_{i;\Sigma_2}M_{\sigma(i);\Sigma_1}^{-1}
X_{\sigma(i);\Sigma_1}.
$$
Let us prove that
the monomial $M(y_1,\dots,y_n)=M_{i;\Sigma_2}M_{\sigma(i);\Sigma_1}^{-1}$
is, in fact, trivial, that is, equals~$1$.

Indeed, assume that there exists a variable, say, $y_1$, that enters 
$M$ with a negative exponent $-k$, $k>0$. Consider the $(n+1)\times
2n$-matrix 
$\widehat{B}$ obtained from $(B_\pr)_{\Sigma_1}$ by adding the 
$(n+1)$th row
$(-d_1b_{1,n+1;\Sigma_1},\dots,-d_nb_{n,n+1;\Sigma_1},0,\dots,0)$, 
where $d_1,\dots,d_n$ are the diagonal entries of $D$.

This matrix, together with the cluster 
$\widehat{\x}=(x_{1;\Sigma_1},\dots, x_{n;\Sigma_1},y_1)$ defines a 
cluster algebra $\A(\widehat{B})$ of rank $n+1$ with $m-1$ stable 
variables $y_2,\dots,y_m$. Put 
$\widehat{\Sigma}_1=( \widehat{\x},\widehat{B})$. The map 
$\ \widehat{\cdot}\ $ is naturally extended to all seeds of
$\A_\pr(B)$ 
so that if $\Sigma,\Sigma'$ are adjacent in direction $k\in [1,n]$ 
then $\widehat{\Sigma},\widehat{\Sigma}'$ are adjacent in the same direction.
Clearly, $\A_\pr(B)$ is the restriction of  $\A(\widehat{B})$ to the 
first $n$ variables. Therefore, one can use only mutations in directions $1,\dots,n$ in
$\A(\widehat{B})$ to get
$$
x_{i;\widehat{\Sigma}_2}=x_{n+1;\widehat{\Sigma}_1}^{-k}
M'(y_2,\dots,y_n)x_{\sigma(i);\widehat{\Sigma}_1},
$$
where $x_{n+1;\Sigma_1}$ is naturally identified with $y_1$. The rest
of the 
proof proceeds exactly as the proof of 
Theorem~\ref{thm:adjclusters}.

\bibliographystyle{amsalpha}

\end{document}